\documentclass[12pt]{article}
\usepackage{graphicx,amsmath,amsfonts,amsthm,verbatim}
\newtheorem{theorem}{Theorem}
\newtheorem{lemma}{Lemma}
\newtheorem{prop}{Proposition}

\theoremstyle{remark}
\newtheorem*{remark}{Remark}

\newcommand{\R}{\mathbb R}
\newcommand{\C}{\mathbb C}
\newcommand\Real{\operatorname{Re}}

\newcommand\Lpsip{L^\psi_+}
\newcommand\Lpsim{L^\psi_-}
\newcommand\Lpsipm{L^\psi_\pm}
\newcommand\RR{\mathbb{R}}
\newcommand\Dom{\operatorname{Dom}}
\newcommand\upto{\uparrow}
\newcommand\downto{\downarrow}
\newcommand\ev{{\operatorname{ev}}}
\newcommand\odd{{\operatorname{odd}}}
\newcommand\qtr{\frac1{4}}
\newcommand\ang[1]{\langle #1 \rangle}
\newcommand\ind{\operatorname{ind}}

\begin{document}

\title{{\bf Eigenvalues of Schr\"odinger operators with potential asymptotically homogeneous of degree $-2$}}

\author{\\
{\bf Andrew Hassell}
\\
\small Department of Mathematics\\
\small The Australian National University\\
\small ACT, Australia
\and
\\
{\bf Simon Marshall\footnote{Address for correspondence: Department 
of Mathematics,
The University of Auckland,
Private Bag 92019,
Auckland New Zealand.
Fax: 64-9-3737457,
email: smar141@ec.auckland.ac.nz }}
\\
\small Department of Mathematics\\
\small The University of Auckland\\
\small New Zealand
\\
}

\maketitle

\normalsize

\begin{abstract}
We strengthen and generalise a result of Kirsch and Simon on the behaviour of the function $N_L(E)$, the number of bound states of the operator $L = \Delta+V$ in $\R^d$ below $-E$. Here $V$ is a bounded potential behaving asymptotically like $P(\omega)r^{-2}$ where $P$ is a function on the sphere. It is well known that the eigenvalues of such an operator are all nonpositive, and accumulate only at $0$.  If the operator $\Delta_{S^{d-1}}+P$ on the sphere has negative eigenvalues $-\mu_1,\ldots,-\mu_n$ less than $-(d-2)^2/4$, we prove that $N_L(E)$  may be estimated as
\[
N_L(E)) = \frac{\log(E^{-1})}{2\pi}\sum_{i=1}^n \sqrt{\mu_i-(d-2)^2/4} +O(1);
\]
thus, in particular, if there are no such negative eigenvalues then $L$ has a finite discrete spectrum. Moreover, under some additional assumptions including that $d=3$ and that there is exactly one eigenvalue $-\mu_1$ less than $-1/4$, with all others $> -1/4$, we show that the negative spectrum is asymptotic to a geometric progression with ratio $\exp(-2\pi/\sqrt{\mu_1 - \qtr})$. 
\vfill

\end{abstract}
\vfill
\newpage

\section{Introduction}
Consider a Schr\"odinger operator of the form $H = \Delta + V$ acting on $L^2(\RR^d)$ for $d \geq 3$, where $\Delta = -\sum_i \partial_i^2$ is the positive Laplacian and $V$ is a multiplication by a potential function $V(x)$ which tends to zero at infinity. 
It is a basic problem in quantum mechanics to determine the location and the nature of the spectrum of $H$. It is well known, and not hard to show, that on $(-\infty, 0)$ the spectrum is discrete. 
It is known that under rather general conditions (decay of $V$ at some prescribed rate at infinity, $O(r^{-1})$ being sufficient), there are no positive eigenvalues of $H$ \cite{Ho}, \cite{FHHH}.  Finally, under fairly general conditions, e.g. $V = O(r^{-1 - \epsilon})$ for some $\epsilon > 0$, the restriction of $H$ to its continuous spectral subspace is unitarily equivalent to $\Delta$ (see e.g. Kato \cite{Kato} or Yafaev \cite{Yafaev}, Theorem 2.4). 
All of these results apply under the assumptions on $V$ made below. The main question, related to the description of the spectrum at this coarse level, left unanswered by the considerations above is whether the point spectrum accumulates at $0$, i.e. whether the point spectrum is finite or infinite. Finding bounds or asymptotics on the negative spectrum has been a major topic of research in Schr\"odinger operators for many decades, and it remains an active field (see e.g. \cite{RoSo}, \cite{GY}). 

Let us denote by $N_H(E)$ the number of eigenvalues of $H$ below $-E<0$ (assuming inf $\sigma_{ess}(H)=0$).  It is well known that the behaviour of $V$ near infinity determines whether $N_H(E)$ is bounded as $E\rightarrow 0$. There is an attractive discussion of this in \cite{RS4}, section XIII.3.  For example, if $V\ge -cr^{-2-\epsilon}$ for some positive $c$ then $N_H(E)$ is bounded, and it is unbounded if $V\le -cr^{-2+\epsilon}$ (here $r=|x|$ is distance to the origin).  In the borderline case where $V$ behaves asymptotically as $-cr^{-2}$, then $N_H(E)$ is bounded if $c<1/4(d-2)^2$ (in fact, the operator $\Delta -(d-2)^2r^{-2}/4$ on $L^2(\RR^d)$ is a positive operator) and unbounded if $c>(d-2)^2/4$, where $d \geq 2$ is the spatial dimension \cite{KS}.  

In this paper we shall analyze the more general case where $V\sim cr^{-2}$ along any ray from the origin, but the constant $c$ depends on  direction, i.e. is   a function on the unit sphere $S^{d-1}$. We may then write $V$ as $r^{-2}(P(\omega)+o(1))$.  Our first main result is

\begin{theorem}\label{main}
Consider the Schr\"odinger operator $L=\Delta+V$ in $\R^d$, where $V$ is a bounded measurable potential equal to $r^{-2} (P(\omega) + t(r, \omega))$ for $r \geq 1$. Denote the eigenvalues of  the operator $\Delta_{S^2}+P$ on $L^2(S^{d-1})$  less than $-\qtr(d-2)^2$ by $-\mu_1,\ldots,-\mu_n$. Assume one of the following two conditions:

(i) The value $-\qtr(d-2)^2$ is not an eigenvalue of  $\Delta_{S^{d-1}}+P$ on $L^2(S^{d-1})$ and the function $t(r, \omega)$ is $O((\log r)^{-(1+ \epsilon)})$ as $r \to \infty$ for some $\epsilon > 0$, or

(ii) $t(r, \omega) = O((\log r)^{-2-\epsilon})$ as $r \to \infty$, for some $\epsilon > 0$.

Then the counting function of $L$ may be estimated as
\begin{equation}
\label{large}
N_L(E) =  \frac{\log(E^{-1})}{2\pi}\sum_{i=1}^n \sqrt{\mu_i-(d-2)^2/4} +O(1)
\end{equation}
In particular, if there are no $-\mu_i\le -\qtr(d-2)^2$ then $L$ has a finite discrete spectrum.

\end{theorem}

\begin{remark} The exponent in condition  (ii) is sharp, in the following sense: there exists a potential $V$ of the form of the form $r^{-2}(-\qtr(d-2)^2 + t(r))$ for $r \geq 1$, with $t = O((\log r)^{-2})$, such that $L = \Delta + V$ has an infinite number of negative eigenvalues, and thus fails to satisfy \eqref{large}. See Section~\ref{examples} for an example. 
\end{remark}

In \cite{KS}, Kirsch and Simon  proved \eqref{large} (with the weaker error estimate\break $o(\log(E^{-1}))$) in the special case that $P$ is a constant function on the sphere. 


To prove our second theorem, we make additional assumptions. We work in three dimensions and assume that $V$ is smooth, that for $r \geq r_0 > 0$ the potential $V$ is exactly equal to $P(\Omega) r^{-2}$, and that the smallest eigenvalue $-\mu_1$ is less than $-1/4$ while all others are strictly greater than $-1/4$.  Under these assumptions we can prove

\begin{theorem}\label{main2} Let $d=3$ and assume that the potential $V$ satisfies the conditions of the previous paragraph. Let
\begin{equation}
\sigma = \exp \Big( -  \frac{2\pi}{\sqrt{ \mu_1 - 1/4}} \Big) .
\label{sigma-defn}\end{equation}
Then there exists $a \in \RR$ such that, with $-E_n$ the $n$th eigenvalue of $\Delta + V$ counted with multiplicity, 
$$
\frac{E_n}{\sigma^n} \to a \ \text{ as } \ n \to \infty. 
$$
Moreover, if $v_n$ is the corresponding normalized eigenfunction, then $v_n$ is essentially supported in an annulus with radii proportional to $\sigma^{-n/2}$, in the following sense: for any $\epsilon > 0$, there exist $C_+, C_- > 0$ so that for all $n \geq 1$, the $L^2$ mass of $v_n$ in the annulus 
$$
C_- \sigma^{-n/2} \leq r \leq C_+ \sigma^{-n/2}
$$
is at least $1-\epsilon$. 
\end{theorem}

The strategy of the proof of Theorem~\ref{main} is as follows. We introduce operators $L_+$ and $L_-$ which lie respectively above and below $L$, and then further modify these to `pseudo-Laplacians' $\Lpsip$ and $\Lpsim$ in order to eliminate difficulties near the origin.  These pseudo-Laplacians $\Lpsipm$ have a domain which is a finite-dimensional perturbation of the domain of $L_\pm$, in the sense that the intersection of the two domains is a closed subspace of finite codimension in each domain. It is a simple consequence of the min-max characterization of eigenvalues that the difference of the counting functions $|N_{\Lpsipm}(E) - N_{L_\pm}(E)|$ is bounded by the maximum of these codimensions. Hence, to estimate $N_L(E)$ up to $O(1)$ it suffices to get asymptotics for $N_{\Lpsip}(E)$ and $N_{\Lpsim}(E)$ which differ by $O(1)$ as $E \to 0$. The operators $\Lpsipm$ are arranged so that we can separate variables by introducing the eigendecomposition on the sphere,  and we thereby  reduce the problem to a very classical problem of estimating the number of zeros of the solution to an ordinary differential equation in the $r$ variable. 
Our pseudo-Laplacians are somewhat analogous to operators used elsewhere in spectral theory to control difficulties caused by  low-lying eigenvalues on a cross-section, as in \cite{LP},  \cite{CdV}, \cite{CZ} for example. 

\

To prove Theorem~\ref{main2}, we construct approximate eigenfunctions $\Phi_n$ for a sequence of approximate eigenvalues $\lambda_n$ asymptotic to $a \sigma^n$, where $\sigma$ is given by \eqref{sigma-defn}. For this to be effective we need the norm of the error term $(H - \lambda_n) \Phi_n$ to be $o(\lambda^{n}) \| \Phi_n \|_2$ as $n \to \infty$. Our approximate eigenfunction $\Phi_n$  is equal to  a zero-mode of $\Delta + V$ inside a ball of radius $\rho_n > r_0$, and an exact solution of the eigenfunction equation outside the ball, where we may use separation of variables, taking advantage of the assumption that $V$ is exactly homogeneous in this region. In addition we need to add smoothing terms so that $\Phi_n$ is $C^1$ across the interface.  The approximate eigenvalues $\lambda_n$ are determined by the criterion that the \emph{principal} terms in the expansions of $\Phi_n$ inside and outside the ball $B(0, \rho_n)$ match in a $C^1$ way at the interface; this automatically makes the $\lambda_n$ asymptotic to a geometric progression (Lemma~\ref{approxev}). 

The result of the construction is that we can show that there is at \emph{least} one eigenvalue in intervals of the form $[ a \sigma^n - o(\sigma^n), a \sigma^n + o(\sigma^n)]$; note that these are non-overlapping intervals for large $n$. By Theorem~\ref{main}, however, there can be at \emph{most} one eigenvalue in all but finitely many of these intervals, and at most finitely many eigenvalues lying outside the union of these intervals. We conclude that in some spectral interval $[-\lambda_0, 0)$ there is exactly one eigenvalue in each of these intervals and no other eigenvalues. It then follows that the $\Phi_n$ are close to true eigenfunctions, and the statement about the $L^2$ mass of the eigenfunctions can be verified by checking it for the $\Phi_n$. 

\

In the final section we consider two examples. The first illustrates the remark after Theorem~\ref{main}. The second example shows some limitations of a heuristic approach of Fefferman and Phong to obtaining estimates for the counting function $N_H(E)$  for an operator $H$. Their approach is to count the number of disjoint images of the unit cube in phase space $\RR^n \times \RR^n$, under canonical transformations, which can be fitted into the part of phase space where the symbol of $H$ is less than $E$. While this heuristic has been shown to give accurate estimates of counting functions (up to constants depending only dimension) in many cases \cite{FP1}, \cite{FP2}, our example shows that no such heuristic involving \emph{connected} regions of phase space will work in the present setting. 

\section{Auxiliary Operators}
Our method is to estimate the operator $L$, and hence $N_L(E)$, from above and below by two operators $L_\pm$ with a simplified potential. We then compare these operators with two further `pseudo-Laplacians' $\Lpsipm$ which allow us to treat the `near region' $r \leq 1$ and the `far region' $r \geq 1$ separately. Let us choose a constant $W$ such that $|V(x)| \leq W$ for all $x \in \RR^d$. Let 
$$t(r) = C (1 + \log r)^{-(1+\epsilon)}, $$
where $\epsilon$ is as in Theorem~\ref{main} and $C$ is large enough so that $|t(r, \omega)| \leq t(r)$ for $r \geq 1$\footnote{Under the assumption (ii), we will take $t(r) = C' (1 + \log r)^{-(2 + \epsilon)}$}.
 We now define the potential operators $V_\pm$ by
\begin{equation}
\label{potential}
V_\pm(x)= |x|^{-2} P(\omega) \pm t(r), \quad |x| \ge 1; \qquad
V_\pm(x) = \pm W, \quad |x| <1
\end{equation}
and introduce the operators $L_\pm$ given by
$$
L_\pm = \Delta + V_\pm
$$
with domain $H^2(\RR^d)$. Clearly we have
\begin{equation}
L_- \leq L \leq L_+.
\label{comparison}\end{equation}

The pseudo-Laplacians $\Lpsipm$ are defined by the same formula as $L_\pm$ but with a different domain. The domain of $\Lpsipm$ consists of functions whose restrictions to the interior and exterior of the unit ball $B^d$ are in $W_2^2(B^d)$ and $W_2^2(\R^d/B^d)$ respectively, together with some conditions at the boundary of the ball which depends on a choice of finite-dimensional subspace $T$ of $L^2(S^{d-1})$. Since the sphere $S^{d-1}_r$ of radius $r$ centred at the origin is canonically identified with the sphere $S^{d-1}$ of radius $1$, by radial projection, we can regard $T$ as a finite dimensional subspace in $L^2(S^{d-1}_r)$ for any positive $r$. Let $\Pi_T$, $\Pi_T^\perp$ denote the orthogonal projections onto $T$, and onto its orthogonal complement, respectively, in $L^2(S^{d-1}_r)$. 

To describe these conditions, first note that $\psi|_{B^d}\in W_2^2(B^d)$ and $\psi|_{\R^d/B^d}\in W_2^2(\R^d/B^d)$ implies, by the Sobolev embedding theorem, that
 the limits $\lim_{r \downto 1} \psi(r, \cdot)$, $\lim_{r \downto 1} \partial_r \psi(r, \cdot)$, $\lim_{r \upto 1} \psi(r, \cdot)$, $\lim_{r \upto 1} \partial_r \psi(r, \cdot)$ are all well-defined in $L^2(S^{d-1})$. Indeed, $\psi(r, \cdot)$ is a $C^{1, 1/2}$ function of $r$ with values in $L^2(S^{d-1})$. The conditions we place on $\psi \in \Dom(\Lpsipm)$ is that 
\begin{equation}\begin{gathered}
\Pi_T ( \lim_{r \downto 1} \psi(r, \cdot) ) = \Pi_T ( \lim_{r \upto 1} \psi(r, \cdot) ) = 0, \\
\Pi_T^\perp (   \psi(r, \cdot) ) \in W^2_2(\RR^3)  \\
\end{gathered}\label{Tcond}\end{equation}

Let us now prove 
\begin{lemma} The operators  $\Lpsipm$ are self-adjoint. 
\end{lemma}

\begin{proof}
From the definition of adjoint, we seek the operator $(\Lpsipm)^*$ with maximal domain such that
\begin{equation}
\label{adj}
\int (\Delta+V_\pm)\psi\overline{\phi} = \int \psi\overline{(\Lpsipm)^*\phi}
\end{equation}
for all $\psi\in \Dom(\Lpsipm)$.  If we restrict to $\psi$ supported inside $B^d$, respectively $\R^d \setminus B^d$, we see that $\phi$ has a generalised Laplacian on these domains so $\phi|_{B^d}\in W_2^2(B^d)$ and $\phi|_{\R^d \setminus B^d}\in W_2^2(\R^d \setminus B^d)$.  We may then apply Stokes' formula to obtain
\begin{multline*}
\int_{\RR^d} \psi\overline{(\Lpsipm)^*\phi} - \int_{\RR^d} \psi\overline{(\Delta+V_\pm)\phi} \ = \ \underset{r \uparrow 1}{\lim}\left( \int_{S^2_r} \frac{\partial\psi}{\partial {\bf n}}\overline{\phi} - \int_{S^2_r} \psi\overline{\frac{\partial\phi}{\partial {\bf n}}} \right) \\
- \underset{r \downarrow 1}{\lim}\left( \int_{S^2_r} \frac{\partial\psi}{\partial {\bf n}}\overline{\phi} - \int_{S^2_r} \psi\overline{\frac{\partial\phi}{\partial {\bf n}}} \right).
\end{multline*}
As usual this means that the left and right hand sides separately vanish. Thus $(\Lpsipm)^* \phi = (\Delta + V_\pm) \phi$. To determine the domain, we consider the boundary terms at $S^{d-1}$ and write $\psi$ and $\phi$ as the sum of their projections into $T$ and into $T^\perp$. The cross terms cancel by orthogonality of $T$ and $T^\perp$, and we can analyze the terms in $T$ and in $T^\perp$ separately. For the part in $T^\perp$, the values of $\psi$ and $\partial_{\bf n} \psi$ agree from $r\leq 1$ and $r \geq 1$, and this forces the values of $\phi$ and  $\partial_{\bf n} \psi$ to agree also. This implies that $\Pi_T^\perp \phi$ is in $W^{2}_2(\RR^d)$.  For the part in $T$, we have $\psi = 0$ on  both sides, but the values of $\partial_{\bf n} \psi$ are independent. This forces $\phi$ to behave likewise. This proves that $\Dom(\Lpsipm)^* = \Dom(\Lpsipm)$. 
\end{proof}

The domain $H^2(\RR^3)$ of $L_\pm$ can be described as those $\psi$ with $\psi|_{D^3}\in W_2^2(D^3)$ and $\psi|_{\R^3/D^3}\in W_2^2(\R^3/D^3)$ satisfying the second condition in \eqref{Tcond}, with the first replaced by
\begin{equation}\begin{gathered}
\Pi_T (  \lim_{r \downto 1} \psi(r, \cdot) ) = \Pi_T ( \lim_{r \upto 1} \psi(r, \cdot) ) , \\
\Pi_T (  \lim_{r \downto 1} \partial_r \psi(r, \cdot) ) = \Pi_T ( \lim_{r \upto 1} \partial_r \psi(r, \cdot) ) 
\end{gathered}\end{equation}
(in other words, the function $\Pi_T(\psi)$ and its normal derivative $\Pi(\partial_r \psi)$ are consistent across the unit sphere). 
This makes it clear that the intersection $\Dom L_\pm \cap \Dom \Lpsipm$ is a closed subspace of codimension $\tau = \dim T$ inside both $\Dom L_\pm$ and inside $\Dom \Lpsipm$ (with respect to the graph norm). It follows from this and from the minimax characterisation of eigenvalues that
\begin{equation}
|N_{L_\pm}(E) - N_{\Lpsipm}(E)| \leq \tau. 
\label{tau}\end{equation}
Combining \eqref{comparison} and \eqref{tau}, we see that to prove \eqref{large} it suffices to prove that both $\Lpsip$ and $\Lpsim$ satisfy the asymptotic on the right hand side of \eqref{large}.

\section{Characterisation of Eigenfunctions of $\Lpsipm$}

We shall determine $N_{\Lpsipm}(E)$ up to a constant by exactly describing all but a finite number of the eigenfunctions of $\Lpsipm$. 

If $\psi$ is such an eigenfunction of $\Lpsipm$ with eigenvalue $-\lambda< 0$, for $r\ge 1$ we may expand $\psi$ in the eigenfunctions $\{H_i\}_{i=1}^{\infty}$ of $\Delta_{S^{d-1}}+P(\omega)$,
\[
\psi(r,\omega)=\sum_{i=1}^\infty X_i(r)H_i(\omega),
\]
and for $r<1$ we may expand in the eigenfunctions of $\{J_i\}_{i=1}^{\infty}$ of $\Delta_{S^{d-1}}\pm W$,
\[
\psi(r,\omega)=\sum_{i=1}^\infty Y_i(r)J_i(\omega).
\]
The functions $\{H_i\}_{i=1}^\infty$ and $\{J_i\}_{i=1}^\infty$ are ordered by eigenvalue, and we recall the the eigenvalues of $\Delta_{S^{d-1}}+P(\omega)$ are $\{-\mu_i\}_{i=1}^\infty$. Let $\nu_i$ denote the $i$th eigenvalue of $\Delta_{S^{d-1}}\pm W$.  Because the potential $V_\pm$ is smooth in $r$ for $r \neq 1$  the $X_i$ and $Y_i$ are smooth, and we may separate variables to obtain the following equations:
\begin{equation}
\label{ODE1}
r^2X_i''+(d-1)rX_i'+(\mu_i \mp t(r)-r^2\lambda)X_i=0, \quad r \geq 1
\end{equation}
\begin{equation}
\label{ODE2}
r^2Y_i''+(d-1)rY_i'-(\nu_i+r^2\lambda)Y_i=0, \quad r \leq 1.
\end{equation}
Because $\psi\in L^2(\R^3)$ we must have $X_i(r)\rightarrow 0$ as $r\rightarrow \infty$.  We claim that if $T$ is chosen to contain enough of the functions $H_i$ and $J_i$, the eigenfunctions of $\Lpsipm$ are relatively simple.  In particular, we wish to prove the following:

\begin{lemma}
There is a $k$ such that if $T=\text{span}\{H_1,\ldots,H_k,J_1\ldots,J_k\}$, all but a finite number of eigenfunctions of $\Lpsipm$ are of the form
\begin{equation}
\psi(r,\omega) = \sum_{i=1}^k X_i(r)H_i(\omega), \quad r\ge 1,
\label{X}\end{equation}
\begin{equation}
\psi(r,\omega)= \sum_{i=1}^k Y_i(r) J_i(\omega), \quad r\le 1.
\label{Y}\end{equation}
\end{lemma}

\begin{proof} 
We begin by recalling that  the nonpositive spectrum of $L^\psi_\pm$ is discrete except for a possible accumulation point at $0$. 
Because we may disregard as many eigenfunctions as we like, we may assume that the eigenvalue $-\lambda$ is restricted to an interval $[-\lambda_0, 0]$, for some $\lambda_0 > 0$. Our strategy is  to show that when $\lambda_0$ is sufficiently small, there exists $k$ such that for  $i \geq k$, the inequalities $Y_i(1)Y_i'(1)\ge 0$ and $X_i(1)X_i'(1) \le 0$ are satisfied, with equality only if $X_i$ or $Y_i$ are identically $0$.

This may be proven for any $k$ with $-\mu_k \geq \sup_{r \geq 1} |t(r)|$.  
If $i > k$, $X_i(1)\ge 0$ and $X_i'(1)>0$, and $X_i'(r)$ is not positive for all subsequent $r$, there must have been some first $r_0$ for which $X_i'(r_0)=0$.  We know $X_i(r_0)>0$, and substituting $r_0$ into (\ref{ODE1}) we obtain
\[
r^2X_i''(r_0)=(-\mu_i \pm t(r)+r^2\lambda)X_i(r_0).
\]
We know $-\mu_i \pm t(r) \geq 0$  and $r^2\lambda>0$, so this implies $X_i''(r_0)>0$ which contradicts the assumption that $X_i'$ was becoming negative for the first time.  Therefore $X_i'$, and so $X_i$, will be positive for all $r$, contradicting the assumption that $X_i(r) \rightarrow 0$.  We may obtain a similar contradiction in the case $X_i\le 0$, $X_i'<0$ so we have $X_i(1)X_i'(1)\le 0$.  In the case $X_i'(1)=0$, $X_i(1)\neq 0$ we may show that $X_i''(1)$ and $X_i(1)$ must have the same sign, so $X_i$ has a positive local maximum (negative local minimum) and we again get a contradiction.  Therefore equality may only occur when $X_i(1)=X_i'(1)=0$.

We now wish to show that $Y_i(1)Y_i'(1)\ge 0$ for all but finitely many $i$.  Because our eigenfunction is continuous at the origin, and each spherical harmonic (except for the constant function, which we may exclude by taking $k\ge 1$) has average value $0$, their coefficients in the expansion of $\psi$ must vanish at zero.

Choose $k$ so that $\nu_k$, the $k$th eigenvalue of $\Delta_{S^{d-1}} \pm W$ is greater than $0$.  If $Y_i(1)\ge0$ and $Y_i'(1)<0$, then $Y$ must have a positive local maximum somewhere on the interval $(0,1)$. However substituting such a point into (\ref{ODE2}) gives us $\nu_i<-r^2\lambda < 0$,  a contradiction.  The case $Y_i(1)\le 0$ and $Y_i'(1)>0$ is similar.  Finally, if $Y_i'(1)=0$ then $Y_i(1)$ and $Y_i''(1)$ will have the same signs, in which case we must again have either a positive local maximum or negative local minimum on $(0,1)$.  Therefore $Y_i(1)Y_i'(1)\ge 0$ with equality only if $Y_i(1)=Y_i'(1)=0$, i.e. $Y_i$ vanishes on $[0,1]$.

To combine $Y_i(1)Y_i'(1)\ge 0$ and $X_i(1)X_i'(1) \le 0$ and obtain a contradiction, we observe that
\begin{equation}
\int_{\partial\Omega} \Pi^\perp_T  \psi\frac{\partial(\Pi^\perp_T \psi)}{\partial {\bf n}}dV
\label{exp}\end{equation}
is the same when the derivative is taken from the inside and outside of the unit sphere. This follows from the fact that $\Pi^\perp_T \psi \in W^2_2(\RR^d)$.
Expanding in eigenfunctions from the inside and from the outside then gives 
\begin{eqnarray*}
0 \ge \sum_{i=k+1}^\infty X_i(1)X_i'(1) & = & \sum_{i=k+1}^\infty Y_i(1)Y_i'(1) \ge 0
\end{eqnarray*}

Therefore $X_i(1)X_i'(1)=0$ and $Y_i(1)Y_i'(1)=0$ for all $i> k$, so only the spherical functions contained in $T$ may be present in any of our remaining eigenfunctions, proving the lemma. 
\end{proof}

\section{A Sturm-Liouville Reformulation}

Lemma $1$ effectively reduces the study of the discrete spectrum of $\Lpsipm$ to the finite number of one dimensional Sturm-Liouville problems \eqref{X}, \eqref{Y}. We may rewrite (\ref{ODE1}) and (\ref{ODE2}) as

\begin{equation}
\label{SL1}
\Gamma_{\mu_i} X_i \equiv X_i''+\frac{d-1}{r}X_i'+\frac{\mu_i+t(r)}{r^2}X_i = \lambda X_i
\end{equation}
and
\begin{equation}
\label{SL2}
\Omega_{\nu_i} Y_i \equiv Y_i''+\frac{d-1}{r}Y_i'-\frac{\nu_i}{r^2}Y_i = \lambda Y_i.
\end{equation}

The $X_i$ have boundary conditions $X_i(1)=0$ (due to the domain condition \eqref{Tcond}) and $\underset{r \rightarrow \infty}{\lim} X_i(r) =0$.  $Y_i$ must satisfy $Y_i(0)=Y_i(1)=0$ for $i>1$ and $Y_1(1)=Y_1'(0)=0$, the Neumann condition at $0$ again arising from smoothness of $\psi$ at the origin.  To analyse the spectrum of $\Gamma$ and $\Omega$ we use the following well-known link between the number of eigenvalues of the operator $\Gamma_{\mu}$ less than $E$ and the number of zeroes of a solution to \eqref{SL1} with $\lambda = E$; see for example \cite{RS4}, Theorem XIII.8.

\begin{theorem}
If $n(\Gamma_{\mu},E)$ is the number of zeros of the solutions of $\Gamma_{\mu}X=E X$, $X(1)=0$ for $E>0$, and $N(\Gamma_{\mu},E)$ the counting function of the operator $\Gamma_{\mu}$ with boundary conditions as above, then $N(\Gamma_{\mu_i},E)=n(\Gamma_{\mu_i},E)$ for  $E > 0$.  The analogous result holds for the $\Omega_{\nu_i}$.
\end{theorem}

As we are considering only a finite number of operators, the $O(1)$ may be taken to be uniform.  It remains to estimate the number of zeros of a solution to $\Gamma_{\mu}X=E X$, $X(1)=0$.  To do so we transform $\Gamma_{\mu}$ by substitution: setting $h(r)=Xr^{(d-2)/2}$ we obtain the equation
\[
r^2h''+rh'+(\mu- \qtr (d-2)^2-r^2E+t(r))h=0,
\]
and setting $g(r)=h(e^r)$,
\begin{equation}
g''+\Big(\mu-\qtr (d-2)^2-E e^{2r}+t(e^r)\Big)g=0.
\label{g}\end{equation}

\emph{Proof under hypothesis (i)}. In this case, by assumption there is no eigenvalue equal to $-\qtr (d-2)^2$ so we can split our analysis of \eqref{g} into the cases $-\mu<-\qtr (d-2)^2$ and $-\mu>-\qtr (d-2)^2$.  In the first case, normalise the coefficient $\mu-\qtr (d-2)^2$ to be 1 by letting $\alpha=\frac{1}{\sqrt{\mu-(d-2)^2/4}}$ and defining $f(r)=g(\alpha x)$. The final form of our equation is then

\begin{equation}
\label{trans1}
f''+\Big(1-\alpha^2E e^{2\alpha r}+t(e^{\alpha r})\Big)f=0
\end{equation}

We may think of this equation as describing paths $(f,f')$ in a phase plane, so that estimating the number of roots becomes a question of estimating the rate at which our point travels around the origin.  Let $\theta$ be the angle between the point $(f,f')$ and the positive $x$ axis.  We have
\begin{eqnarray*}
\theta & = & \tan^{-1}\frac{f'}{f}\\
\implies \frac{d\theta}{dr} & = & \frac{1}{1+(f'/f)^2}.\frac{f''f-(f')^2}{f^2}\\
& = & \frac{f''f-(f')^2}{f^2+(f')^2}\\
& = & -\frac{(1-\alpha^2E e^{2\alpha r}+t(e^{\alpha r}))f^2+(f')^2}{f^2+(f')^2}\\
& = & -1+\frac{(\alpha^2E e^{2\alpha r}-t(e^{\alpha r}))f^2}{f^2+(f')^2}\\
\implies |\theta'+1| & \le & \alpha^2E e^{2\alpha r}+t(e^{\alpha r})
\end{eqnarray*}
Therefore for small $r$, the point $(f,f')$ will travel clockwise around the origin at a nearly constant rate.  The error term $\alpha^2E e^{2\alpha r}+t(e^{\alpha r})$ becomes very great for large $r$, but this is not a problem because as soon as $1-\alpha^2E e^{2\alpha r}+t(e^{\alpha r})$ becomes negative the character of the ODE \eqref{trans1} changes.  In particular, if $r_0$ is such that $1-\alpha^2E e^{2\alpha r}+t(e^{\alpha r})<0$ for $r>r_0$, then by \eqref{trans1} $f$ cannot have a positive local maximum (or negative local minimum) with $r>r_0$ and so can only have one root greater than $r_0$.  Therefore to estimate the number of times $(f,f')$ winds around the origin we only need to consider the behaviour on the interval $[1,r_0]$ where $1+\|t\|_{\infty} \ge \alpha^2E e^{2\alpha r}$ holds.  For all $E$
\begin{equation}
\int_0^{r_0} \alpha^2E e^{2\alpha r}dr  \le  \frac{1}{2}\alpha Ee^{2\alpha r_0} \le  \frac{1+\|t\|_{\infty}}{2\alpha},
\end{equation}
so the contribution to the error from the $\alpha^2E e^{2\alpha r}$ term is bounded.  The form of $t(r)$ implies that $\int_0^\infty t(e^{\alpha r})dr < \infty$ giving us
\[
\theta(r_0)-\theta(0)=-r_0+O(1),
\]
so the number of roots of $f$ is $r_0/\pi+O(1)$.  Calculating $r_0$ in terms of $E$ (up to a O(1) term) gives $r_0 = \log(1/E)/2\alpha +O(1)$. Therefore we have the number of roots equal to $\frac{1}{2\pi}\sqrt{\mu-1/4}\ln(1/E)+O(1)$. (We also remark that, since there are only finitely many eigenvalues below $-d(d-2)^2/4$, there are only a finite number of values of $\alpha$ and so the $O(1)$ estimates are uniform in $\alpha$.)

In the case $-1/4<-\mu$, we have $1/4-\mu>0$ and so may rewrite our equation for $g$ as
\[
g''-(1/4-\mu+E e^{2r}-t(e^r))g=0.
\]
Because $t=o(1)$, there will be an $r_0$ independent of $E$ such that $(1/4-\mu+E e^{2r}-t(e^r))>0$ for $r>r_0$ and $E$ small, and after which $g$ can have at most one more root as discussed above.  Also, the number of roots that $g$ may have in the interval $r<r_0$ is bounded for small $E$ by the continuity of solutions of ODEs in their coefficients in compact domains.  Therefore when $-1/4<-\mu$, the number of roots is bounded as $E\rightarrow 0$.

A similar analysis of the operators $\Omega_i$ is much easier, because they act on functions with compact domains and so solutions to $\Omega_i Y=E Y$, $Y(1)=0$ are uniformly continuous in $E$.  Therefore the number of roots of solutions is bounded as $E\rightarrow 0$, so these problems only have a finite number of positive eigenvalues.  

  We therefore have
\[
N_L(E)) = \frac{\log(E^{-1})}{2\pi}\sum_{i=1}^n \sqrt{4\mu_i-1} +O(1)
\]
as required.  The reverse estimate on the potential, and the counting function, is performed in exactly the same way.

\emph{Proof under hypothesis (ii)} In this case, we have to analyze equation \eqref{g} also in the case that $-\mu = -1/4$, but now we take that $t(r) = C (\log r)^{-(2 + \epsilon)}$. Let $T(r) = t(e^r)$, so that $|T(r)| \leq C r^{-(2 + \epsilon)}$.  We have to show that the number of zeroes of the solution of \eqref{g} is $O(1)$ as $E \to 0$. 

Since the $-E e^{2r}$ term is negative, by the Sturm comparison theorem it suffices to show that the number of zeroes of the solution of
\begin{equation}
g'' + T g = 0
\label{Tg}\end{equation}
is finite. Also, by the Sturm comparison theorem, we may assume that $T \geq 0$ (otherwise, we can consider $\max(T, 0)$). 

To do this, let $q$ be a root of $g$, and assume (by scaling $g$) that $g'(q) = 1$. (We cannot have $g'(q) = 0$ unless is the trivial solution.) To show that $g$ has no root larger than $q$, it suffices to show that 
$$
\int_q^\infty -g''(r) \, dr \leq 1,
$$
since that will show that $g$ is increasing on $[q, \infty)$. We have
$$
\int_q^\infty -g''(r) \, dr \leq \int_q^\infty T g \, dr \leq C \int_q^\infty r^{-(2 + \epsilon)} (r-q) \, dr \leq \frac{C}{\epsilon} q^{-\epsilon}.
$$
This is indeed less than $1$ for $q$ sufficiently large. 
Thus for $r_0$ sufficiently large, depending on $C$ and $\epsilon$, we see that $g$ can have at most one root on $[r_0, \infty)$. The number of roots in the interval $[1, r_0]$ is uniformly bounded for $E \in [0, E_0]$ so this gives a uniform bound on the whole real line. This completes the proof of the theorem under hypothesis (ii).


\section{Preliminaries to the proof of Theorem~\ref{main2}}

Our method of proof is the technique of approximate eigenfunctions, described in the following well-known lemma:

\begin{lemma}\label{approx}
Let $A$ be a self-adjoint operator on a Hilbert space $H$.
 Suppose there exists
$u \neq 0 \in H$ with 
\begin{equation}
 \| (A - \lambda) u \| \leq \epsilon \| u \|.
\label{cond}\end{equation}
Then there exists a point of spec$(A)$ in the interval $[\lambda -
\epsilon, \lambda + \epsilon]$.

Now let $\mu > \epsilon$, and suppose that the interval $[\lambda -
\mu, \lambda + \mu]$ contains only one point of spec$(A)$ which is an
eigenvalue of multiplicity $1$. Let $v$ be a normalized eigenfunction
for this eigenvalue. Assume also that $\langle u, v \rangle \geq 0$
(if not, then replace $v$ by $\alpha v$ where $\alpha \in \C$,
$|\alpha| = 1$ is chosen so that $\langle u, \alpha v \rangle \geq
0$). Then if $u$ satisfies \eqref{cond},
$$
\big\| \frac{u}{\| u \|} - v \big\| \leq \frac{6\epsilon}{\mu}.
$$
\end{lemma}

\begin{proof} If there is no spectrum of $A$ in the interval
$[\lambda - \epsilon, \lambda + \epsilon]$ then, since the spectrum
is closed, there is $\epsilon' > \epsilon$ such that there is no
spectrum of $A$ in the interval $[\lambda - \epsilon', \lambda +
\epsilon']$. By the spectral mapping theorem, then, $\| (A -
\lambda)^{-1} \| \leq 1/\epsilon'$. This contradicts \eqref{cond}.
Hence spec$(A) \cap [\lambda - \epsilon, \lambda + \epsilon]$ must be
nonempty.

To prove the second part of the lemma, we write $u = cv + w$ where $c
= \langle u, v \rangle > 0$  and $w$ is orthogonal to $v$. Then
$$
\| (A - \lambda) u \| \geq \| (A - \lambda) w \| - \| (A - \lambda) cv \| \\
\geq \mu \| w \| - c \epsilon
$$
which implies that
$$
\mu \| w \| \leq (1 + c) \epsilon \implies \| w \| \leq \frac{2\epsilon}{\mu}.
$$
Now
$$
\| u - v \| \leq \|  w \| + \| v - cv \| \leq \frac{2\epsilon}{\mu} + 1 - c.
$$
Finally
$$
\| w \|^2 = 1 - c^2 \implies 1 - c = \frac{\| w \|^2}{1 + c} \leq
\big( \frac{2\epsilon}{\mu} \big)^2 \leq \frac{4\epsilon}{\mu}
$$
so
$$
\| u - v \| \leq  \frac{6\epsilon}{\mu}.
$$
\end{proof}

We therefore wish to find a sequence of approximate eigenvalues $-\lambda_n$, tending to zero,  and  functions $\Phi = \Phi_n$ such that $\|(L+\lambda_n)\Phi_n\|/\| \Phi_n \|$ is small.  Because $V$ is unperturbed for $r \ge r_0$, we may separate variables to construct exact solutions of $L\Phi = \lambda_n\Phi$ in this region.  However, these exact solutions are very poor approximations where $V$ is perturbed.  One way around this would be to use a given eigenfunction of of $L$ in this region, however if we choose an eigenfunction $\Psi$ with eigenvalue $-\nu < 0$ for this purpose we have $(L+\lambda_n)\Phi = (\lambda_n-\nu)\Psi$ in the region where $\Psi$ is used.  This gives us an error with a $|\lambda_n-\nu|$ term in it which remains rougly constant as $\lambda_n \rightarrow 0$, which is bad from the point of view of applying Lemma~\ref{approx}.  A better choice for $\Psi$ is a well behaved zero-mode of $L$, because now $(L+\lambda_n)\Phi = \lambda_n\Psi$ and we have an error behaving as $\lambda_n$.  The existence of such a function is guaranteed by the following result:

\begin{prop}\label{zero}
Let the operator $\Delta+V$ in $\R^3$ be as described above. Let $-\mu_1, \mu_1 > 1/4$ be the smallest eigenvalue, and $J_1(\omega)$ the corresponding eigenfunction, of the spherical operator $\Delta_{S^2}+P$.  
Denote by  $N_\epsilon$ denote the null space of $\Delta + V$ acting on $\ang{x}^{-\epsilon} L^2(\R^3)$, where $\ang{x} = (1 + |x|^2)^{1/2}$.  
Then for small positive $\epsilon$ the dimension of $N_{-1-\epsilon} / N_{-1+\epsilon}$ is one dimensional. This quotient space is spanned by a function $\Psi$ satisfying 
$$
(\Delta+V)\Psi = 0
$$
and having the asymptotic expansion 
$$
\Psi = r^{-1/2}\cos(\sqrt{\mu_1-1/4}\ln(r)+c)J_1(\omega) + O(r^{-1/2-\delta}), \quad r \to \infty
$$
for some $c \in \R$ and some $\delta > 0$. 
\end{prop}

\begin{proof} These results are implied by the relative index theorem in chapter 6 of \cite{Mel}. Let us write $r$ for a function that is $\geq 1$ everywhere on $\R^3$ and equal to $|x|$ for $|x| \geq \rho_0$. We first note that we can write
$$
\Delta + V = r^{-5/2} \Big( -(r \partial_r)^2 + \Delta_{S^2} + P + \frac1{4} \Big) r^{1/2}
$$
for $r > \rho_0$. Therefore the null space of $\Delta + V$ on $\ang{x}^{-1 \pm \epsilon} L^2(\R^3)$ is equal to the null space of an operator $Q$ on $r^{\pm \epsilon} L^2(\RR_+ \times S^2; \mu_b)$ where $Q$ takes the form $ -(r \partial_r)^2 + \Delta_{S^2} + P + \frac1{4}$ near $r = \infty$ and $\mu_b$ is a smooth measure taking the form $d\omega dr/r$ for $r \geq \rho_0$ (where $d\omega$ is the standard measure on $S^2$). The operator $Q$ is an elliptic b-operator in the sense of \cite{Mel}. By the relative index theorem of \cite{Mel}, $Q$ is Fredholm as a map\footnote{Here $H^2_b$ is the b-Sobolev space of order two, given by the functions whose derivatives with respect to angular derivatives $\partial_\omega$ and with respect to $r \partial_r$ up to second order are square integrable.}
\begin{equation}
Q : r^{- \epsilon} H^2_b(\RR_+ \times S^2; \mu_b) \to r^{- \epsilon} L^2(\RR_+ \times S^2; \mu_b)
\label{QQ}\end{equation}
for all $\epsilon \neq 0$ in a neighbourhood of $0$,  and, with $N^Q_\epsilon$ denoting the null space and $\ind_\epsilon Q$ the index of $Q$ acting as in \eqref{QQ}, for small $\epsilon > 0$,
$$
\ind_{-\epsilon} Q - \ind_\epsilon Q = 2.
$$
(This number is equal to the number of indicial roots, i.e. complex numbers $\tau$ for which $\Delta_{S^2} + P + 1/4 + \tau^2$ is not invertible, with imaginary part in the interval $[-\epsilon, \epsilon]$. Here we have $\tau = \pm \sqrt{\mu_1 - 1/4}$ as the only such indicial roots, provided $\epsilon$ is sufficiently small.)
Since $Q$ is self-adjoint on $L^2(\RR_+ \times S^2; \mu_b)$ and by elliptic regularity we have 
$$
\ind_\epsilon Q = \dim N^Q_\epsilon - \dim N^Q_{-\epsilon} .
$$
Therefore $\ind_\epsilon Q = - \ind_{-\epsilon} Q$, and we have
for sufficiently small $\epsilon > 0$
$$
\ind_{-\epsilon} Q   = 1. 
$$
This implies that $\dim ( N_{-1-\epsilon} / N_{-1+\epsilon} ) = 1$. If we take any function $\Psi$ in $N_{-1-\epsilon} \setminus N_{-1+\epsilon}$, then the regularity results of \cite{Mel}, e.g. Proposition 5.21, show that $\Psi$ has a conormal asymptotic expansion as $r \to \infty$, i.e. an expansion in powers and powers of logarithms of $r$ with coefficients that are eigenfunctions of the operator on the boundary. This implies the asymptotic property stated in the theorem.  
\end{proof} 


Our approximate eigenfunctions, constructed in the following section, will be built out of the function $\Psi$ from Proposition~\ref{zero}, close to the origin, and an exact solution to the eigenfunction equation in the far region, where we can use separation of variables. 
If the appropriate spherical eigenfunctions are $\{J_i\}_{i=1}^\infty$ with eigenvalues $\{-\mu_i\}_{i=1}^\infty$, we may expand $F$ as

\[
F(r,\omega) = \sum_{i=1}^\infty f_i(r)J_i(\omega),
\]

where for $r>r_0$ the components $f_i$ obey the differential equation
\begin{equation}
\label{ODE}
r^2f_i''+2rf_i'+(\mu_i-r^2\lambda)f_i=0.
\end{equation}

Letting $f_i = r^{-1/2} \tilde f_i$, the function $\tilde f_i$ satisfies the modified Bessel equation
\begin{equation}
\label{ODE'}
r^2 \tilde f_i''+r \tilde f_i'+(\mu_i - \qtr-r^2\lambda)\tilde f_i=0.
\end{equation}
There is a one-dimensional space of solutions which are exponentially decreasing as $r \to \infty$, given by multiples of the MacDonald functions $\tilde f_i = K_{\sqrt{1/4 - \mu_i}}(\sqrt{\lambda} r)$ (see \cite{AS}, section 9.6). It is a standard fact that these functions $K_\nu(r)$ have no zeroes for $r > 0$ (\cite{AS}, p377).  

Let us define the exponents $\alpha_i$, $\beta_i$ by
\begin{equation}
\alpha_i = \frac{-1 - \sqrt{1 - 4\mu_i}}{2}, \quad \beta_i = \frac{-1 + \sqrt{1 - 4\mu_i}}{2}.
\label{ab}\end{equation}
Then the functions $r^{\alpha_i}$, $r^{\beta_i}$ solve \eqref{ODE} with $\lambda = 0$. Let 
\begin{equation}
X_i(r) = r^{\alpha_i} \int_{-\infty}^\infty (p^2+1)^{\alpha_i} e^{-ipr} dp.
\end{equation}
Then $X_i$ is a multiple of $r^{-1/2} K_{\sqrt{1/4 - \mu_i}}(r)$ and solves \eqref{ODE}. 
For each $i\ge 2$,  let 
\begin{equation}
\label{lambda}
X_i^{\lambda}(r)=X_i(r\sqrt{\lambda});
\end{equation}
this is an exponentially decaying solution 
of (\ref{ODE}) for arbitrary $\lambda$, with the explicit expression
\begin{equation}
\label{fourier}
X_i^\lambda = r^{\alpha_i} \lambda^{-\alpha_i/2-1/2} \int_{-\infty}^\infty (p^2+\lambda)^{\alpha_i} e^{-ipr} dp . 
\end{equation}

We now obtain some estimates on the $X_i^\lambda$ needed in the proof of Theorem~\ref{main2}. The identity (\ref{fourier}) lets us estimate the value of $X_i^{\lambda}$ for $\alpha_i < -1/2$, i.e. $i \geq 2$, by
\begin{equation}
\label{Xvalue}
X_i^\lambda (r) = (r \sqrt{\lambda})^{\alpha_i} \Big[ \int_{-\infty}^\infty (p^2 + 1)^{\alpha_i} \, dp \Big] (1+o_{\alpha_i}(1)), \quad r \sqrt{\lambda} \to 0.
\end{equation}
Here the $o$ term is uniform as $i \to \infty$, but blows up as $\alpha_i 
\to -1/2$ (which corresponds to $\mu_i \to -1/4$); however, since there are only a finite number of $\alpha_i$ in this range, for a given potential, we may take this estimate to be uniform in $i \geq 2$.  This implies that
\begin{equation}
X_i^\lambda (r) \geq C (r \sqrt{\lambda})^{\alpha_i} |2 \alpha_i + 1|^{-1} 
\text{ for sufficiently small } r\sqrt{\lambda} . 
\label{Xvalue2}\end{equation}

For $\alpha_i < -3/2$ we may use \eqref{fourier} to prove
\begin{equation}
\frac{(X_i^{\lambda})'(r)}{X_i^{\lambda}(r)}=\frac{\alpha_i}{r}(1+O((r \sqrt{\lambda})^2)), \quad r \sqrt{\lambda} \to 0,
\label{Xiless}\end{equation}
where the $O$ term is uniform as $i \to \infty$.  

On the other hand, series solutions at $r=0$ show that if $-1/2 < \alpha_i \leq -3/2$ and $\alpha_i$ is real,
\begin{equation}
\label{ratioest}
\frac{(X_i^{\lambda})'(r)}{X_i^{\lambda}(r)}=\frac{\alpha_i}{r}(1+O_{\alpha_i}(\ln(r\sqrt{\lambda})(r \sqrt{\lambda})^{-1-2\alpha_i})), \quad r \sqrt{\lambda} \to 0. 
\end{equation}
Here the $O_{\alpha_i}$ is not uniform as $\alpha_i \to -1/2$, but as above, since there are only a finite number of $\alpha_i$ in the range $[-3/2, -1/2)$  for a given potential, we may take this estimate to be uniform in $i$.  If we define 
\begin{equation}
K =  \min\big( 2, \min\limits_{i \geq 2} (-\frac1{2} -  \alpha_i) \big) > 0
\label{Kdefn}
\end{equation}
then we may combine the two estimates \eqref{Xiless} and \eqref{ratioest} in the form
\begin{equation}
\label{ratioest1}
\frac{(X_i^{\lambda})'(r)}{X_i^{\lambda}(r)}=\frac{\alpha_i}{r}(1+O((r \sqrt{\lambda})^{K})), \quad i \geq 2,  \quad r \sqrt{\lambda} \to 0.
\end{equation}

In the case $i=1$, there is a constant $d$ such that
\begin{equation}
\label{X1}
X_1^{\lambda}(r) = r^{-1/2}\lambda^{-1/4}\left( \cos(\sqrt{\mu_1-1/4}\ln(r \sqrt{\lambda})+d)+O(r \sqrt{\lambda}) \right)
\end{equation}
and
\begin{multline}
\label{X1'}
(X_1^{\lambda})'(r) = r^{-3/2}\lambda^{-1/4}\Big( -1/2\cos(\sqrt{\mu_1-1/4}\ln(r \sqrt{\lambda})+d) \\  -\sqrt{\mu_1-1/4}\sin(\sqrt{\mu_1-1/4}\ln(r \sqrt{\lambda})+d)+O(r \sqrt{\lambda}) \Big),
\end{multline}
where both functions $O(r \sqrt{\lambda})$ are smooth away from $0$.  

We also need estimates on $L^2$ norms. Using \eqref{fourier}   we may bound the $L^2$ norm of $X_i^\lambda J_i$ restricted to $\R^3\backslash B(0, \rho)$ as follows:
\begin{equation}\begin{aligned}
\| X_i^\lambda J_i & \|^2_{L^2(\RR^3 \setminus B(0, \rho))}  =  \int_\rho^\infty (X_i(r \sqrt{\lambda}))^2 r^2dr\\
& =  \int_\rho^\infty (r \sqrt{\lambda})^{2\alpha_i} \left( \int_{-\infty}^\infty (p^2+1)^{\alpha_i} e^{-ipr \sqrt{\lambda}} dp \right)^2 r^2dr\\
& \le  \rho^{2\alpha_i+2} \lambda^{\alpha_i-1/2} \int_\rho^\infty \left( \int_{-\infty}^\infty (p^2+1)^{\alpha_i} e^{-ipr \sqrt{\lambda}} dp \right)^2 d(r \sqrt{\lambda})\\
& \le  \rho^{2\alpha_i+2} \lambda^{\alpha_i-1/2} \left\| \int_{-\infty}^\infty (p^2+1)^{\alpha_i} e^{-ipr} dp \right\|^2\\
& =  \rho^{2\alpha_i+2} \lambda^{\alpha_i-1/2} \| (p^2+1)^{\alpha_i} \|^2\\
& \le  O_{\alpha_i}(1) \rho^{2\alpha_i+2} \lambda^{\alpha_i-1/2} \label{zzz}
\end{aligned}\end{equation}
Here the $O(1)$ is uniform as $\alpha_i \to -\infty$, but not as $\alpha_i \to -1/2$; again, as above,  for a given potential, we may take this estimate to be uniform in $i$.


\section{Approximate eigenfunctions}

Our approximate eigenfunctions, for suitable values of the eigenvalue $-\lambda$, $\lambda > 0$, will be built out of three components. These are the zero mode $\Psi$ from Proposition~\ref{zero}, close to the origin; an exact solution to the eigenfunction equation far from the origin; and smoothing terms supported in an intermediate region, which ensures that the two pieces fit together so as to lie in the domain of $L$. To define this, let $h(t)$ be a function which is supported on the interval $[1,2]$, equal to $t-1$ on $[1,1.5]$, and smooth on $(1,2]$. Hence $h$ is smooth except at $t=1$ where it has a jump in the first derivative. 
For the $n$th approximate eigenvalue, $-\lambda_n$, we will choose a radius $\rho_n$, growing as a negative power of $\lambda_n$, and define
\begin{equation}
\Phi = \Phi_n = 
\begin{cases}
\Psi, \quad \hskip 155pt r \leq \rho_n, \\
\sum_{i=1}^\infty \Big( \phi_i X_i^{\lambda_n}(r) + \chi_i h(r/\rho_n) \Big)  J_i(\omega) 
\quad r \geq \rho_n.
\end{cases}
\label{approx-efn}\end{equation}
Here $\phi_i$ is a coefficient chosen to ensure the continuity of $\Phi$ at $r = \rho_n$, and the coefficient $\chi_i$ is chosen so that the radial derivative of $\Phi$ is continuous at $r = \rho_n$. 

We may perform separation of variables on $\Psi$ in the region $r \geq r_0$ to obtain an expansion
\[
\Psi(r,\omega) = \sum_{i=1}^\infty Y_i(r)J_i(\omega).
\]
The functions $Y_i$ must satisfy (\ref{ODE}) for $\lambda=0$ and $r>r_0$, which has general solution $Y_i(r) = ar^{\alpha_i}+br^{\beta_i}$ where $\alpha_i, \beta_i$ are defined in \eqref{ab}.  For $i\ge 2$, $\alpha_i$ and $\beta_i$ are real and $\alpha_i<-1/2<\beta_i$; since $\Psi=O(r^{-1/2})$, we have 
\begin{equation}
Y_i(r) =\psi_i r^{\alpha_i}, \quad r \geq r_0
\label{Yi}\end{equation}
for some coefficients $\psi_i$.

We choose $\lambda_n$ so that the smoothing term for $i=1$ is not needed, i.e. so that we can choose $\chi_1 = 0$. We therefore need to choose $\lambda$ so that 
\begin{equation}
\label{condition}
(X_1^{\lambda})'(\rho) Y_1(\rho) = Y_1'(\rho)X_1^{\lambda}(\rho). 
\end{equation}

\begin{lemma}\label{approxev} Let 
\begin{equation}
\label{xi}
\xi_n = \exp\left( \frac{-2\pi n-C}{\sqrt{\mu_1-\qtr}} \right)
\end{equation}
for some $C \in \R$ and let  $\rho_n$ be a sequence of real numbers $\geq r_0$ with $\rho_n \sqrt{\xi_n} \to 0$. Then there exists a constant $C$ and a sequence $(\lambda_n)$, satisfying (\ref{condition}) at $\rho = \rho_n$ which are asymptotic to $\xi_n$ in the sense that 
\begin{equation}
\frac{\lambda_n}{\xi_n} = 1 + O(\rho_n \sqrt{\xi_n}) \text{ as } n \to \infty.
\label{asympt}\end{equation}
\end{lemma}

\begin{remark} All the $O(\cdot)$  and $o(\cdot)$ estimates in this section are uniform as $i \to \infty$, $\lambda \to 0$ and $\rho \to \infty$ provided $\rho \sqrt{\lambda} \to 0$. 
\end{remark}

\begin{proof}
Since we know that 
$$Y_1(r) = r^{-1/2} \cos( \sqrt{\mu_1 - 1/4} \log r + c), \quad r \geq r_0,
$$
 we have
\begin{equation*}\begin{aligned}
& Y_1'(r)  \\  =  &- r^{-\frac{3}{2}} \Big(  \frac1{2} \cos( \sqrt{\mu_1 - 1/4} \log r + c) + \sqrt{\mu_1 - 1/4} \sin( \sqrt{\mu_1 - 1/4} \log r + c) \\
 =  &- r^{-\frac{3}{2}} \mu_1 \cos( \sqrt{\mu_1 - 1/4} \log r + c - \theta')
\end{aligned}\end{equation*}
where $\cos \theta' = (2 \sqrt{\mu_1})^{-1}$. Similarly we have 
$$
X_1^\lambda(r) = r^{-1/2} \lambda^{-1/4} \Big( \cos( \sqrt{\mu_1 - 1/4} \log (r \sqrt{\lambda}) + d) + O(r \sqrt{\lambda}) \Big)
$$
and 
$$\begin{gathered}
\frac{d}{dr} X_1^\lambda(r) 
= - r^{-3/2} \lambda^{-1/4}  \mu_1 \Big( \cos( \sqrt{\mu_1 - 1/4} \log (r \sqrt{\lambda}) + d - \theta') + O(r \sqrt{\lambda}) \Big) 
\end{gathered}$$
as $r \sqrt{\lambda} \rightarrow 0$.  We want to solve 
$$
Y_1(r) \frac{d}{dr} X_1^\lambda(r) = Y_1'(r)  X_1^\lambda(r)
$$
which is equivalent to
$$\begin{gathered}
\cos( \sqrt{\mu_1 - 1/4} \log r + c)  \Big( \cos( \sqrt{\mu_1 - 1/4} \log (r \sqrt{\lambda}) + d - \theta') + O(r \sqrt{\lambda}) \Big) \\
 = \cos( \sqrt{\mu_1 - 1/4} \log r + c - \theta') \Big( \cos( \sqrt{\mu_1 - 1/4} \log (r \sqrt{\lambda}) + d) + O(r \sqrt{\lambda}) \Big)
\end{gathered}$$
at $r = \rho_n$. Using trigonometric identities this reduces to
$$
\sin \theta'  \sin(\sqrt{\mu_1 - 1/4} \log \sqrt{\lambda} + d - c) = O(\rho_n \sqrt{\lambda}).
$$
The $\sin \theta'$ can be neglected here since this is a positive quantity depending only on $\mu_1$. We note that the $O(\rho_n \sqrt{\lambda})$ term on the right hand side is a continuous function of 
$\lambda$. Therefore, by the intermediate value theorem we get a sequence of solutions of the form
$$
\frac1{2} \log \lambda_n = \frac{c - d - \pi n}{\sqrt{\mu_1 - 1/4}} + O(\rho_n \sqrt{\lambda_n})
$$
for $n$ sufficiently large, which satisfies \eqref{asympt} with $C = 2(d - c)$. 
\end{proof}

We now analyse the matching conditions required at\footnote{For ease of notation we will drop the subscript $n$ from $\lambda_n$ and $\rho_n$ from now on.} $r = \rho$. In order to have continuity of $\Phi$ at $r = \rho$ we choose, for $i \geq 2$, 
\begin{equation}
  \phi_i    = \frac{Y_i(\rho)}{X_i^{\lambda}(\rho)} 
  = \frac{Y_i(\rho)}{X_i(\rho \sqrt{\lambda})}.
\label{phii}\end{equation}
and in order to have continuity of $\partial_r \Phi$ at $r = \rho$ we choose 
$$
\chi_i = \rho \big( Y_i'(\rho)-\phi_i(X_i^{\lambda})'(\rho) \big).
$$
(We remark that $X_i$ has no zeroes for $i \geq 2$, so the denominator in \eqref{phii} is never zero.) Using \eqref{phii}, \eqref{Yi} and \eqref{ratioest1} we find that
\begin{equation}\begin{gathered}
\chi_i = \psi_i \alpha_i \rho^{\alpha_i} - \psi_i \rho^{\alpha_i} \frac{ \rho (X_i^\lambda)'(\rho)}{X_i^\lambda(\rho)} \\
= \psi_i \rho^{\alpha_i} \Big( \alpha_i - \frac{ \rho (X_i^\lambda)'(\rho)}{X_i^\lambda(\rho)} \Big) = \psi_i \rho^{\alpha_i} O(|\alpha_i|(\rho \sqrt{\lambda})^{K}).
\end{gathered}\label{chiiest}\end{equation}
where $K$ is defined in \eqref{Kdefn}.

\begin{lemma}
The expression \eqref{approx-efn}  converges in the graph norm on the domain of $L$.
\end{lemma}

\begin{proof}
Since $\Psi$ is smooth, it is enough to show that $\Phi$ is in the Sobolev space $W^2_2(\RR^2 \setminus B(0, r_0))$. In this set we 
have $\Phi = \sum_i s_i$, where $s_i$ is given by
$$
s_i = \begin{cases} \psi_i Y_i(r) J_i(\omega), \quad r \leq \rho \\
\big( \phi_i X_i^\lambda(r) + \chi_i h(r/\rho) \big) J_i(\omega), \quad r \geq r_0.
\end{cases}
$$
Due to our matching conditions, each $s_i$ has two derivatives in $L^\infty$, and $s_i$ is rapidly decreasing together with all derivatives as $r \to \infty$. Therefore, each term $s_i$ is in the domain of $L$. We have to show that the sum converges in the graph norm. Since $\Psi \in W^2_2(B(0, \rho))$ and $L \Psi = 0$, it is sufficient to show convergence of the infinite sums 
\begin{equation}
\label{sum1}
\sum_{i=2}^\infty \int_\rho^\infty (\phi_i X_i^\lambda + \chi_i h_\rho)^2r^2dr
\end{equation}
and
\begin{equation}
\label{sum2}
\sum_{i=2}^\infty \| (L + \lambda)(\phi_i X_i^\lambda + \chi_i h_\rho) J_i \|^2 = \sum_{i=1}^\infty \int_\rho^{2\rho} \chi_i^2 ( h_\rho''+\frac{2}{r} h_\rho'+\left( \frac{\mu_i}{r^2}-\lambda \right) h_\rho)^2r^2dr
\end{equation}
where we write $h_\rho(r) = h(r/\rho)$. 
We may further reduce (\ref{sum1}) to showing the convergence of
\begin{equation}
\sum_{i=2}^\infty \int_\rho^\infty (\phi_i X_i^\lambda)^2r^2dr \  \ \hbox{and} \ \ \sum_{i=2}^\infty \int_\rho^\infty \chi_i^2 h_\rho^2r^2dr.
\label{twocomp}\end{equation}
Similarly, \eqref{sum2} reduces to showing convergence of 
\begin{equation}\begin{gathered}
\sum_{i=2}^\infty \int_\rho^{2\rho}  \big( \chi_i (h_\rho'' + \frac{h_\rho'}{r}) \big)^2 r^2 \, dr, \qquad
\sum_{i=2}^\infty \int_\rho^{2\rho}
 \big( \frac{\chi_i \mu_i h_\rho}{r^2} \big)^2 r^2 \, dr  \\
\text{ and } \sum_{i=2}^\infty \int_\rho^{2\rho} (\lambda \chi_i  h_\rho)^2  r^2 \, dr
\end{gathered}\label{threecomp}\end{equation}
(recalling that $h_\rho$ is supported in $[\rho, 2\rho]$). 
We may obtain useful bounds on the sums in terms of $\| \Psi_{2r_0} \|$ 
where $\Psi_{2r_0}$ is the restriction of $\Psi$ to the disk of radius $2r_0$, using the estimate
\begin{equation}
\label{L2est}
\int_{r_0}^{2r_0} r^{2\alpha_i}r^2dr \ge C2^{\alpha_i}r_0^{2\alpha_i}
\end{equation}
for some constant $C$ independent of $i$.  Equation (\ref{L2est}) gives us
\begin{equation}
\label{dirich}
\infty > \| \Psi_{2r_0} \|^2 \ge \sum_{i=1}^\infty \int_{r_0}^{2r_0} \psi_i^2r^{2\alpha_i}r^2dr \ge C\sum_{i=1}^\infty 2^{\alpha_i}\psi_i^2 r_0^{2\alpha_i},
\end{equation}
This gives us in particular the convergence of the right hand side, and this allows us to show the convergence of the second component of (\ref{sum1}).

If $n$ is such that $\alpha_n < -3/2$ then we can then estimate, using \eqref{chiiest},
\begin{equation}\begin{aligned}
\sum_{i=n+1}^\infty \int_\rho^\infty (\chi_i h_\rho)^2 r^2dr & \le  \sum_{i=n+1}^\infty \psi_i^2 \rho^{2\alpha_i} \alpha_i^2 O((\rho\sqrt{\lambda})^{2K}) \int_\rho^\infty h_\rho^2 r^2 \, dr \\
& =  \sum_{i=n+1}^\infty \psi_i^2 \rho^{2\alpha_i+3} \alpha_i^2 O((\rho\sqrt{\lambda})^{2K})\\
& \le  \sum_{i=n+1}^\infty \psi_i^2 (2r_0)^{2\alpha_i+3} \alpha_i^2 O((\rho\sqrt{\lambda})^{2K}),  \\
& \le  O((\rho\sqrt{\lambda})^{2K}) \sup (2^{\alpha_i} \alpha_i^2) \sum_{i=n+1}^\infty \psi_i^2 2^{\alpha_i} r_0^{2\alpha_i}\\
& \le  O((\rho\sqrt{\lambda})^{2K}) \sup (2^{\alpha_i} \alpha_i^2) \| \Psi_{2r_0} \|^2 \text{ by } \eqref{dirich} \\
& = O((\rho\sqrt{\lambda})^{2K})
\end{aligned}\label{bound1}\end{equation}
which shows convergence of the second term in \eqref{twocomp}. 

A similar estimate  proves convergence for the first and third terms in \eqref{threecomp}. As for the second term, we need to take into account that $|\mu_i| \to \infty$. Using $\sup_i 2^{\alpha_i} \alpha_i^2 \mu_i^2 < \infty$ (which follows from Weyl asymptotics for the $\mu_i$), we can apply the same argument to show convergence of the second term of \eqref{threecomp}. 
 
 It remains to show the convergence of the first term of (\ref{twocomp}). We have 
\begin{equation}\begin{aligned}
\sum_{i=n+1}^\infty  \int_\rho^\infty & (\phi_i X_i^\lambda)^2 r^2dr  \le  O(1) \sum_{i=n+1}^\infty \phi_i^2 \rho^{2\alpha_i+2} \lambda^{\alpha_i-1/2} \text{ by } \eqref{zzz} \\
& =  O(1) \sum_{i=n+1}^\infty \Big( \frac{\psi_i \rho^{\alpha_i} }{ X_i^\lambda(\rho) } \Big)^2 \rho^{2\alpha_i+2} \lambda^{\alpha_i-1/2} \text{ by } \eqref{phii} \\
& =  O(1) \sum_{i=n+1}^\infty (\rho \sqrt{\lambda})^{-2\alpha_i} \alpha_i^2 \psi_i^2 \rho^{4\alpha_i+2} \lambda^{\alpha_i-1/2} \text{ by } \eqref{Xvalue2} \\
& =  O(1) \sum_{i=n+1}^\infty \psi_i^2 \alpha_i^2 \lambda^{-1/2}  \rho^{2\alpha_i+2} \\
& =  O(\lambda^{-1/2})
\end{aligned}\label{bound2}\end{equation}
using the same reasoning as in \eqref{bound1}. 
\end{proof}

In order to apply Lemma~\ref{approx} we need accurate bounds on the norms of $\Phi$ and $(L + \lambda)\Phi$. 

\begin{lemma} Assume that $\rho = \lambda^{-1/2 + \delta/2}$, where $\delta > 0$ is sufficiently small. Let $\Phi'$ be $\Phi - \phi_1 X_1^\lambda J_1$ restricted to $\RR^3 \setminus B(0, \rho)$. Then for some $\epsilon > 0$ (depending on $\delta$) we have estimates 
\begin{equation}
\| \Phi \|_2 \geq c \lambda^{-1/2},
\label{Phifrombelow}\end{equation}
\begin{equation}
\| \Phi' \|_2 = O(\lambda^{-1/2 + \epsilon})
\label{Phi'}\end{equation}
and
\begin{equation}
\label{total}
\|(L+\lambda)\Phi\|_2=
O(\lambda^{1/2 + \epsilon}). 
\end{equation}
\end{lemma}

\begin{proof} To obtain a lower bound on $\| \Phi \|_2$ it is enough to estimate $$\| \phi_1  X_1^\lambda J_1 \|_{L^2(\RR^3 \setminus B(0, \rho)}$$ in view of orthogonality of the $J_i$ on the sphere. By scaling, we compute
\begin{equation}\begin{gathered}
\int_\rho^\infty (X_1^\lambda(r))^2 r^2 \, dr = \lambda^{-3/2} \int_{\rho \sqrt{\lambda}}^\infty X_1^2 (r) r^2 \, dr \\
 \geq c \lambda^{-3/2} \text{ since } \rho \sqrt{\lambda} \to 0.
 \end{gathered}\end{equation}
On the other hand, we may estimate $\phi_1$ by combining the equations 
\begin{eqnarray*}
\phi_1 X_1^\lambda(\rho)=Y_1(\rho),\\
\phi_1 (X_1^\lambda)'(\rho)=Y_1'(\rho)
\end{eqnarray*}
to get
\[
\phi_1^2(X_1^\lambda(\rho)^2+\rho^2(X_1^\lambda)'(\rho)^2)=Y_1(\rho)^2+\rho^2 Y_1'(\rho)^2.
\]
From (\ref{X1}) and (\ref{X1'}) we see that the left hand side is bounded above and below by nonzero constant multiples of $\phi_1^2\rho^{-1}\lambda^{-1/2}$ and the right hand side is similarly estimated by $\rho^{-1}$ (note that  $\psi_1 \neq 0$ so the lower bound is nonzero), so $\phi_1 \sim \lambda^{1/4}$.  Therefore the norm of $\phi_1X_1^{\lambda}J_1$ behaves as $\lambda^{-1/2}$.

To prove \eqref{Phi'}, we note that \eqref{bound1} and \eqref{bound2} show that if we take the sum of $s_i$  over $i \geq n+1$ instead of $i \geq 2$ then we get a bound  $O(\lambda^{-1/2})$ for the square of the $L^2$ norm. To bound the terms with $2 \leq i \leq n$ we can adapt these estimates. In the second line of \eqref{bound1}, if we remove a factor of $\rho^2$ then the rest of the argument follows, since we have $2\alpha_i + 1 < 0$ for $i \geq 2$. Therefore we can end up with $O(\rho^{2})$ for the square of the norm of each $s_i$, which is $O(\lambda^{-1 + \delta})$. In \eqref{bound2} we can reach the second last line, and then we observe that for $2 \leq i \leq n$, the factor  $\rho^{2\alpha_i + 2}$ is bounded by $\rho^{2\alpha_2 + 2}$, giving an overall estimate of $O(\lambda^{-1/2} \rho^{2\alpha_2 + 2})$ for the square of the norm of $s_i$. This is $O(\lambda^{-1 + 2\epsilon})$ for $\epsilon = (-2\alpha_2 - 1)/4 > 0$  which is also of the claimed form. 

To prove \eqref{total} we must estimate the sums in \eqref{threecomp}. To do this we use the computation \eqref{bound1}, modifying it as needed. For the first and second terms in \eqref{threecomp} we have an extra factor of $\rho^{-4}$, as compared to the computation in \eqref{bound1}, so the sum of these terms (for $i \geq n+1$) is $O(\rho^{-4})$. In the case of the third term we have an extra factor of $\lambda^2$, so this term contributes $O(\lambda^2)$. 

Of course we have to bound the terms where $i \leq n$. Since there are only a finite number of these terms we can ignore the $\alpha_i$, $\mu_i$  and $\psi_i$. Thus, for the first and second type of term in \eqref{threecomp},  for a single $i \leq n$, we get a bound (from the second line of \eqref{bound1}, recalling we have an extra $\rho^{-4}$)
$O(\rho^{2\alpha_i - 1} (\rho \sqrt{\lambda})^{2K}) = O(\rho^{2\alpha_2 - 1} (\rho \sqrt{\lambda})^{2K})$. For the third term we get $O(\lambda^2 \rho^{2\alpha_2 + 3} (\rho \sqrt{\lambda})^{2K})$. 
In summary, we obtain 
$$
\| (L + \lambda) \Phi \|^2 = O(\rho^{-4}) + O(\lambda^2) + O(\rho^{2\alpha_2 - 1} (\rho \sqrt{\lambda})^{2K}) + O(\lambda^2 \rho^{2\alpha_2 + 3} (\rho \sqrt{\lambda})^{2K}).
$$
Since $\rho \sqrt{\lambda} \to 0$ and $K > 0$ this is
$
O(\rho^{-4}) + O(\rho^{2\alpha_2 - 1} ) .
$
Now we choose $\rho\sqrt{\lambda} = \sqrt{\lambda}^{\delta}$, i.e. $\rho = \lambda^{(\delta - 1/2)/2}$ where $0 < \delta < 1/2$ is chosen so small that 
$$
\rho^{-4} = O(\lambda^{1 + 2\epsilon}) \text{ and } \rho^{2\alpha_2 - 1} = O(\lambda^{1 + 2\epsilon}).
$$
This requires $\delta < \min(1/4, -\alpha_2(1-2\alpha_2)^{-1})$. With this choice of $\rho$ as a function of $\lambda$ we have $\| (L + \lambda) \Phi \|^2 = O(\lambda^{1 + 2\epsilon})$ and hence \eqref{total}. 
\end{proof}


\section{Proof of Theorem~\ref{main2}}
Applying Lemma~\ref{approx} to $\Phi^{\lambda_n}$, we see that  estimates \eqref{Phifrombelow} and \eqref{total} imply there is a point of spec$(L)$ in the interval $[-\lambda_n-O(\lambda_n^{1+\epsilon}),-\lambda_n+O(\lambda_n^{1+\epsilon})]$.  Combined with Lemma~\ref{approxev}, this shows that for all $n \ge n_0$ \begin{equation}
\label{interval}
\text{spec} \, L \cap [-\xi_n-O(\xi_n^{1+\epsilon}), -\xi_n+O(\xi_n^{1+\epsilon})] \text{ is nonempty}
\end{equation}
where $\xi_n$ and $n_0$ are as in the lemma.  As discussed in the introduction, we may combine this with Theorem~\ref{main} to show that in some spectral interval $[-\lambda_0,0)$ there is exactly one eigenvalue in each of these intervals and no other eigenvalues.  Consequently if $-E_n$ is the $n^{\rm th}$ eigenvalue of $L$, counted with multiplicity, then
\[
\underset{n \rightarrow \infty}{\lim} \frac{E_n}{\xi_n}
\]
converges.  This is equivalent to saying that if
\[
\sigma = \exp \Big( -  \frac{2\pi}{\sqrt{ \mu_1 - 1/4}} \Big)
\]
then the limit
\begin{equation}
\label{thm2}
\underset{n \rightarrow \infty}{\lim} \frac{E_n}{\sigma^n}
\end{equation}
converges, which proves the first part of Theorem~\ref{main2}.

We now wish to apply the second part of Lemma~\ref{approx} to show that if $v_n$ is the normalised eigenfunction corresponding to $-E_n \in [-\lambda_0,0)$ then $\|\tilde\Phi^{\lambda_n}-v_n\|$ is small. Here $\tilde \Phi^{\lambda_n} = \Phi^{\lambda_n} / \| \Phi^{\lambda_n}  \|_2$ is the $L^2$- normalized function.  We have just seen that for $n$ and $C$ sufficiently large there is exactly one eigenvalue in the interval $[-\xi_n - C \xi^{1 + \epsilon}, -\xi_n + C \xi^{1 + \epsilon}]$ and no other eigenvalues in the interval $[-\lambda_0, 0)$. Equivalently, there are integers $n_0$,  $n_1$ such that
$$
E_{n} \in [\xi_{n+n_1} - C \xi_{n+n_1}^{1 + \epsilon}, \xi_{n+n_1} + C \xi_{n+n_1}^{1 + \epsilon}] \text{ for all } n \geq n_0.
$$
Lemma~\ref{approx} gives 
\begin{equation}
\label{eigenfn}
\|\tilde\Phi^{\lambda_{n + n_1}}-v_n\| = O(\xi_n^{\epsilon}) \text{ for  all } n \geq n_0.
\end{equation}
By \eqref{Phifrombelow} and \eqref{Phi'}, the $L^2$ mass of $\tilde\Phi^{\lambda_n}$ is essentially supported outside $B(0, \rho_n)$ by the first spherical eigenfunction $J_1$.  Moreover, since $|\Psi| \leq c r^{-1/2}$, the $L^2$ mass of $\tilde\Phi^{\lambda_n}$ in $B(0, \rho_n)$ is $O(\rho_n \sqrt{\lambda_n}) \to 0$. Finally, because $X_1^\lambda$ is a scaled version of the fixed function $X_1$, and since $\lambda_n/\sigma_n$ converges by \eqref{asympt},  there exist $C_+$ and $C_-$ such that $1-\epsilon$ of the mass of $X_1^{\lambda_n} J_1$ is supported in the annulus
$$
C_- \sigma^{-n/2} \leq r \leq C_+ \sigma^{-n/2}.
$$
We conclude that there are $C_+$ and $C_-$ such that $1-\epsilon$ of the mass of $\tilde\Phi^{\lambda_{n + n_1}}$ is supported in the annulus
$$
C_- \sigma^{-n/2} \leq r \leq C_+ \sigma^{-n/2}.
$$
Combined with (\ref{eigenfn}), this proves the second part of the theorem.


\section{Examples}\label{examples}

We now give an example illustrating the remark following Theorem~\ref{main}. That is, we show  that there is a $t(r)$ with $t(e^r) = O(r^{-2})$, such that the Schr\"odinger operator with potential equal to $r^{-2}(-\qtr  (d-2)^2 + t(r))$  has an infinite number of eigenvalues. 

Following the reasoning of the previous few paragraphs, it is enough to exhibit a function $T(r)$ which is $O(r^{-2})$ such that the ODE \eqref{Tg} has an infinite number of zeros on $[1, \infty)$; then the number of  zeros of \eqref{g} must grow without bound as $E \to 0$ since, on each compact set, the solution will converge uniformly as $E \to 0$. 

We choose a $C^2$ function $g(t)$ on the interval $[1,4]$ so that $g(1) = 1$, $g(2) = -1$ and $g(4) = 1$, such that $g$ is strictly monotone on the intervals $(1,2)$ and $(2,4)$, and so that $g''$ vanishes on a neighbourhood of the point in $(1,2)$ and the point in $(2,4)$ where $g(t) = 0$. Finally, we assume that $g'(1) = 4 g'(4)$ and $g''(1) = 16 g''(4)$. It is clear that such a function exists. We now define $g$ on $[1, \infty)$ by requiring it to be multiplicatively periodic with period ratio $4$, i.e., that $g(4t) = g(t)$. The function $g(t)$ is now $C^2$ on $[1, \infty)$. The ratio
$$
T(t) = -\frac{g''(t)}{g(t)}
$$
is also a $C^2$ function since we required $g''$ to vanish near where $g = 0$. By the periodicity assumption, $T(t) = O(t^{-2})$. Clearly $g$ solves the ODE \eqref{Tg} and has an infinite number of zeroes. 

\

Our next example is related to the work of Fefferman-Phong. We work in dimension 3 for convenience, though the construction works in all dimensions $d \geq 3$.  Let $(\theta, \phi)$ be the usual coordinates on $S^2$, in which the upper hemisphere is $\{ \theta \leq \pi/2 \}$. 
We consider a potential function $P(\theta)$ on the upper hemisphere of $S^2$ which is equal to $-1/3$ for $\theta < \pi/2 - 2\epsilon$, equal to zero for $\pi/2 - \epsilon < \theta < \pi/2$ and monotone in between. Let $P_\ev$ denote the even continuation of $P$ to the sphere, and $P_\odd$ the odd continuation. We assume that $\epsilon$ is sufficiently small (less than $0.01$, say). 

\begin{lemma} 

(i) The lowest eigenvalue of the operator $\Delta + P_\ev$ on $L^2(S^2)$ is less than $-1/4$. 

 (ii) The lowest eigenvalue of the operator $\Delta + P_\odd$ on $L^2(S^2)$ is greater than $-1/4$. 
 \end{lemma}
 
 \begin{proof} To prove (i), we note that the quadratic form
 $$
 Q_\ev(f) = \int_{S^2} |\nabla f|^2 + P_\ev |f|^2 \, d\omega
 $$
 takes a value $-1/3 + O(\epsilon)$ on the $L^2$-normalized constant function,
 which is $< -1/4$ for $\epsilon$ sufficiently small. 
 
 To prove (ii), we take an arbitrary $f \in L^2(S^2)$ with $\| f \|_2 = 1$, and write it $f = f_\ev + f_\odd$ in terms of its even and odd parts with respect to reflection in the equator $\theta = \pi/2$. Then, due to the symmetries of the sphere, 
$$
  \int_{S^2} |\nabla f|^2  \, d\omega =   \int_{S^2} |\nabla f_\ev|^2  \, d\omega +   \int_{S^2} |\nabla f_\odd|^2  \, d\omega.
$$
  On the other hand, $f_\odd$ is in the domain of the  Dirichlet Laplacian on the upper hemisphere. The lowest eigenvalue of the Dirichlet Laplacian on the upper hemisphere is $2$ (this follows since eigenfunctions of the Dirichlet Laplacian on the upper hemisphere are the restrictions of odd spherical harmonics, and the smallest eigenvalue of an odd spherical harmonic is $2$). Hence 
$$
 \int_{S^2} |\nabla f_\odd|^2  \, d\omega \geq 2  \int_{S^2} | f_\odd|^2  \, d\omega
 $$
 and so 
 $$
  \int_{S^2} |\nabla f|^2  \, d\omega \geq 2 \int_{S^2} | f_\odd|^2  \, d\omega. 
$$
Thus
$$\begin{gathered}
Q_\odd(f) \equiv \int_{S^2} |\nabla f|^2 + P_\odd |f|^2 \, d\omega \\
\geq 2 \| f_\odd \|_2^2 +  \int_{S^2} P_\odd ( |f_\ev|^2 + |f_\odd|^2 + 2 \Real f_\ev \overline{f_\odd} ) \, d\omega.
\end{gathered}$$
Since $P_\odd$ is odd and $|f_\ev|^2$, $|f_\odd|^2$ are even, these terms vanish and we are left with 
$$\begin{gathered}
Q(f) \geq 2 \| f_\odd \|_2^2 -  2 \int_{S^2}  |P_\odd |  |f_\ev||f_\odd|  \, d\omega \\
\geq 2\| f_\odd \|_2^2 - \frac{2}{3} \| f_\ev \|_2 \| f_\odd \|_2 \\
\geq  2\| f_\odd \|_2^2 - \frac1{3} \big( \frac1{6} \| f_\ev \|_2^2  + 6\| f_\odd \|_2^2 \big) \\
\geq -\frac1{18}  \| f_\ev \|_2^2 \geq -\frac1{18}  \| f \|_2^2 
\end{gathered}$$
which shows that the lowest eigenvalue of $\Delta + P_\odd$ is no smaller than $-1/18$. 
\end{proof}

According to Theorem~\ref{main}, then, if we choose potentials $V_\odd$ and $V_\ev$ which are equal to $r^{-2} P_\odd$, respectively $r^{-2} P_\ev$, for $r \geq 1$, and are equal to $+1$ for $r < 1$, then the first operator has a finite number of eigenvalues, while the second has an infinite number. On the other hand, it is clear that the region of phase space where $\sigma(L_\ev) < -E$ is the disjoint union of two identical copies of the region where $\sigma(L_\odd) < -E$. So there are an infinite number of disjoint copies of the unit cube, under canonical transformations, in the second region if and only if there are in the first region. Moreover, this holds regardless of the precise conditions we place on these maps. Hence a Fefferman-Phong-type heuristic is not effective in estimating the number of eigenvalues of both these operators.

\end{document}